\documentclass[12pt,leqno]{amsart}

\usepackage{amsmath,amssymb,latexsym,psfrag}
\usepackage{graphicx} 
\usepackage{amsthm,amscd,amsfonts}

\begin{document}

\newtheorem{thm}{Theorem}[section]
\newtheorem{prop}[thm]{Proposition}
\newtheorem{lem}[thm]{Lemma}
\newtheorem{cor}[thm]{Corollary}
\newtheorem{conj}[thm]{Conjecture}
\newtheorem{ddef}[thm]{Definition}
\newtheorem{ex}[thm]{Example}
\newtheorem{rem}[thm]{Remark}
\newtheorem{notation}[thm]{Notation}

\numberwithin{equation}{section}

\newcommand{\bthm}{\begin{thm}}
\newcommand{\ethm}{\end{thm}}
\newcommand{\blem}{\begin{lem}}
\newcommand{\elem}{\end{lem}}
\newcommand{\bcor}{\begin{cor}}
\newcommand{\ecor}{\end{cor}}
\newcommand{\bprop}{\begin{prop}}
\newcommand{\eprop}{\end{prop}}
\newcommand{\bproof}{\begin{proof}}
\newcommand{\eproof}{\end{proof}}
\newcommand{\bddef}{\begin{ddef}}
\newcommand{\eddef}{\end{ddef}}

\newcommand{\beq}{\begin{equation}} \newcommand{\eeq}{\end{equation}}
\newcommand{\beqs}{\begin{equation*}} \newcommand{\eeqs}{\end{equation*}}
\newcommand{\beqarr}{\begin{eqnarray}} \newcommand{\eeqarr}{\end{eqnarray}}
\newcommand{\beqarrs}{\begin{eqnarray*}} \newcommand{\eeqarrs}{\end{eqnarray*}}
\newcommand{\barr}{\begin{array}} \newcommand{\earr}{\end{array}}
\newcommand{\btab}{\begin{tabular}} \newcommand{\etab}{\end{tabular}}

\newcommand{\bit}{\begin{itemize}} \newcommand{\eit}{\end{itemize}}
\newcommand{\ben}{\begin{enumerate}} \newcommand{\een}{\end{enumerate}}
\newcommand{\bce}{\begin{center}} \newcommand{\ece}{\end{center}}

\newcommand{\defeq}{\stackrel{\rm def}{=}}
\newcommand{\bd}{\partial}
\newcommand{\op}{\mathrm{op}}
\newcommand{\hz}{\widehat{0}}
\newcommand{\ho}{\widehat{1}}

\newcommand{\wh}{\widehat}
\newcommand{\wt}{\widetilde}
\newcommand{\cov}{\lhd}
\newcommand{\ccov}{\rhd}

%%some of the usual math blackboard caps

\newcommand{\N}{\mathbb{N}}
\newcommand{\Z}{\mathbb{Z}}
\newcommand{\C}{\mathbb{C}}
\newcommand{\R}{\mathbb{R}}
\renewcommand{\P}{\mathbb{P}}
\newcommand{\Q}{\mathbb{Q}}
\newcommand{\F}{\mathbb{F}}

%% some script caps

\renewcommand{\AA}{\mathcal{A}}
\newcommand{\BB}{\mathcal{B}}
\newcommand{\CC}{\mathcal{C}}
\newcommand{\DD}{\mathcal{D}}
\newcommand{\EE}{\mathcal{E}}
\newcommand{\FF}{\mathcal{F}}
\newcommand{\II}{\mathcal{I}}
\newcommand{\LL}{\mathcal{L}}
\newcommand{\QQ}{\mathcal{Q}}
\renewcommand{\SS}{\mathcal{S}}

%% Greek

\def\al{\alpha}
\def\be{\beta}
\def\de{\delta}
\def\De{\Delta}
\def\ga{\gamma}
\def\Ga{\Gamma}
\def\la{\lambda}
\def\La{\Lambda}
\def\om{\omega}
\def\Om{\Omega}
\def\ze{\zeta}
\def\vphi{\varphi}
\def\veps{\varepsilon}
\def\si{\sigma}
\def\Si{\Sigma}

\def\N{\mathbb{N}}
\def\P{\mathbb{P}}
\def\Z{\mathbb{Z}}
\def\Q{\mathbb{Q}}
\def\R{\mathbb{R}}
\def\C{\mathbb{C}}
\def\L{\mathbb{L}}

\def\AA{{\mathcal{A}}}
\def\BB{{\mathcal{B}}}
\def\CC{{\mathcal{C}}}
\def\FF{{\mathcal{F}}}
\def\LL{{\mathcal{L}}}
\def\NN{{\mathcal{N}}}
\def\HH{{\mathcal H}}
\def\RR{{\mathcal R}}
\def\MM{{\mathcal M}}

\def\bb{\mathbf{b}}
\def\cc{\mathbf{c}}
\def\ff{\mathbf{f}}
\def\bg{\mathbf{g}}
\def\hh{\mathbf{h}}
\def\kk{\mathbf{k}}

\newcommand{\st}{\,:\,} 
\newcommand{\sbseq}{\subseteq}
\newcommand{\spseq}{\supseteq}
\newcommand{\larr}{\leftarrow}
\newcommand{\rarr}{\rightarrow}
\newcommand{\Larr}{\Leftarrow}
\newcommand{\Rarr}{\Rightarrow}
\newcommand{\lrarr}{\leftrightarrow}
\newcommand{\Lrarr}{\Leftrightarrow}

\def\il{\int\limits}
\def\sbs{\subset}
\def\sbseq{\subseteq}
\def\wh{\widehat}
\def\wt{\widetilde}
\def\oli{\overline}
\def\uli{\underline}
\def\langle{\left<}
\def\rangle{\right>}
\def\Lraw{\Longrightarrow}
\def\lraw{\longrightarrow}
\def\Llaw{\Longleftarrow}
\def\llaw{\longleftarrow}
\def\Llraw{\Longleftrightarrow}
\def\llraw{\longleftrightarrow}

\def\wtx{\underset{\displaystyle{\widetilde{}}}{x}}
\def\wth{\underset{\displaystyle{\widetilde{}}}{h_0}}

\def\({\left(}
\def\){\right)}
\def\no={\,{\,|\!\!\!\!\!=\,\,}}
\def\wt{\widetilde}

\def\rank{\text\rm{rank}}
\def\circm{\circmega}
\def\Lv{\left\Vert}
\def\Rv{\right\Vert}
\def\lan{\langle}
\def\ran{\rangle}
\def\wh{\widehat}
\def\no={\,{\,|\!\!\!\!\!=\,\,}}
\def\sbseq{\subseteq}
\def\circli{\circverline}
\def\aff{\text{\rm{aff}}}
\def\conv{\text{\rm{conv}}}
\def\sgn{\circperatorname{sgn}}

\def\lk{\mathrm{link}}

\def\CatTop{\sf Top}
\def\id{\rm id}
\def\colimit{\circperatorname{\sf colim}}
\def\hocolim{\circperatorname{\sf hocolim}}
\def\susp{\circperatorname{susp}}
\def\sd{\circperatorname{sd}}
\def\rk{\mathrm{rk}}
\def\sbseq{\subseteq}
\def\sbs{\subset}
\def\spseq{\supseteq}
\def\sps{\supset}
\def\ssm{\smallsetminus}

\newcommand{\bca}{\begin{cases}} \newcommand{\eca}{\end{cases}}

\newcommand{\OO}{\mathcal{O}}
\newcommand{\LA}{L_{\AA}}
\newcommand{\FA}{F_{\AA}}
\newcommand{\CA}{C_{\AA}}
\newcommand{\MA}{M_{\AA}}
\newcommand{\ep}{\varepsilon}

\newcommand{\mob}{M\"obius function}
\newcommand{\supp}{\mathrm{supp}}
\newcommand{\card}{\mathrm{card}}

\title[Note: random-to-front shuffles on trees]{Note: random-to-front shuffles on trees}
\author[Anders Bj\"orner]{Anders Bj\"orner}\footnote{
%\lowercase{
Research supported by the Knut and Alice Wallenberg Foundation,
grant KAW.2005.0098.}
\address{Royal Institute of Technology, Department of Mathematics,
SE-100 44 Stockholm, Sweden }
\address{  Institut Mittag-Leffler, Aurav\"agen 17,  SE-182 60 Djursholm, Sweden }
\email{bjorner@math.kth.se, bjorner@mittag-leffler.se}
%\thanks{%{\bf Acknowledgement.}}
%\subjclass{}
%\keywords{}

\newcommand{\spa}{\mathrm{span}}
\newcommand{\lang}{\langle\,}
\newcommand{\rang}{\,\rangle}

\begin{abstract} 
A Markov chain is considered whose states are orderings of an underlying fixed tree
and whose transitions are local 
``random-to-front'' reorderings, driven by a probability 
distribution on subsets of
the leaves. The eigenvalues of the transition matrix are determined
using Brown's theory of random walk on semigroups.
\end{abstract}

\maketitle

\section{Introduction}

The random-to-front shuffle of a linear list 
(known in the card-game model also as ``inverse riffle shuffle'')
is a well-known and much studied
finite-state Markov chain. Its states are the linear orderings of an underlying finite set,
and a step of the chain results from selecting a subset (often a singleton)
and moving it to the front of the current list 
in the induced order. See e.g. \cite{BHR,BrDi,FH} and the references given there. 
In this note we consider a slight generalization, namely
to shuffles on trees.

Consider a fixed rooted tree $T$ whose leaves $L$ are all at the same depth.
The following shows a such a tree of depth $3$.
\vspace{.2cm}
\begin{center}
\psfrag{X}{\large $Y$}
\psfrag{XY}{\large $X\circ Y$}
\psfrag{Y}{\large $X$}
\psfrag{x1}{\large {$y_1$}}
\psfrag{x2}{\large $y_2$}
\psfrag{x3}{\large $y_3$}
\psfrag{x4}{\large $y_4$}
\psfrag{x5}{\large {$y_5$}}
\psfrag{x6}{\large $y_6$}
\psfrag{x7}{\large {$y_7$}}
\psfrag{y1}{\large {$x_1 =x_5 =x_7$}}
\psfrag{y2}{\large $x_2 =x_3 =x_4 =x_6$}
\psfrag{R}{\large $\R$}
\resizebox{!}{40mm}{\includegraphics{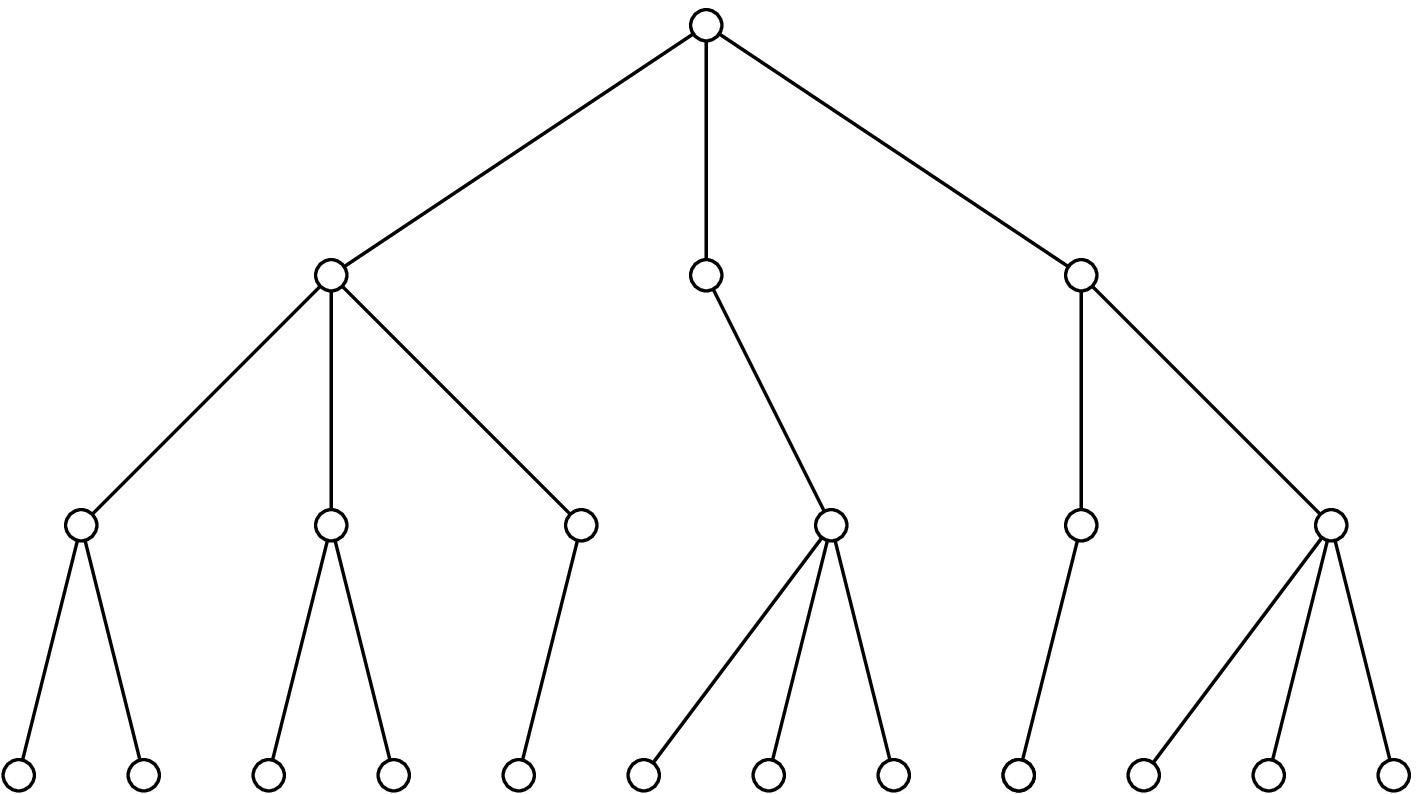}}
%{}{Archimedes}
\end{center}
%\nopagebreak\vspace{.2cm} 
\centerline{{\bf Figure 1.}}
\vspace{3mm}
Suppose that at each inner node 
(i.e., node that is not a leaf)
a total ordering of its children is given.
For instance, it can be the left-to-right ordering given by a planar drawing of the tree,
such as in Figure 1. Now, a subset $E$ of the set of leaves $L$ is
chosen with some probability. Then the ordering is rearranged locally at each inner node
%(the tree is redrawn)
so that the children having some descendant in $E$ come first, and otherwise the induced order
is preserved. The process is illustrated in Figures 3 and 4.

\newpage

In this note the eigenvalues of the transition matrix of this Markov chain are determined.
This is a straight-forward  application of  Brown's theory of random
walks on semigroups \cite{Bro1}. 

\newcommand{\dep}{\text{depth}}

Note that if $\dep(T)=1$ the Markov chain we describe 
amounts  to the classical linear random-to-front shuffle.
% known in the card-game model also as ``reverse riffle-shuffle''. 
For $\dep(T)>1$ we perform such a linear
shuffle locally at each inner node, in each case moving the set of $E$-related nodes
to the front.
\vspace{.2cm}
\begin{center}
\psfrag{ro}{{rooms}}
\psfrag{sh}{{shelves}}
\psfrag{bo}{{books}}
\psfrag{E}{\large E}
\resizebox{!}{50mm}{\includegraphics{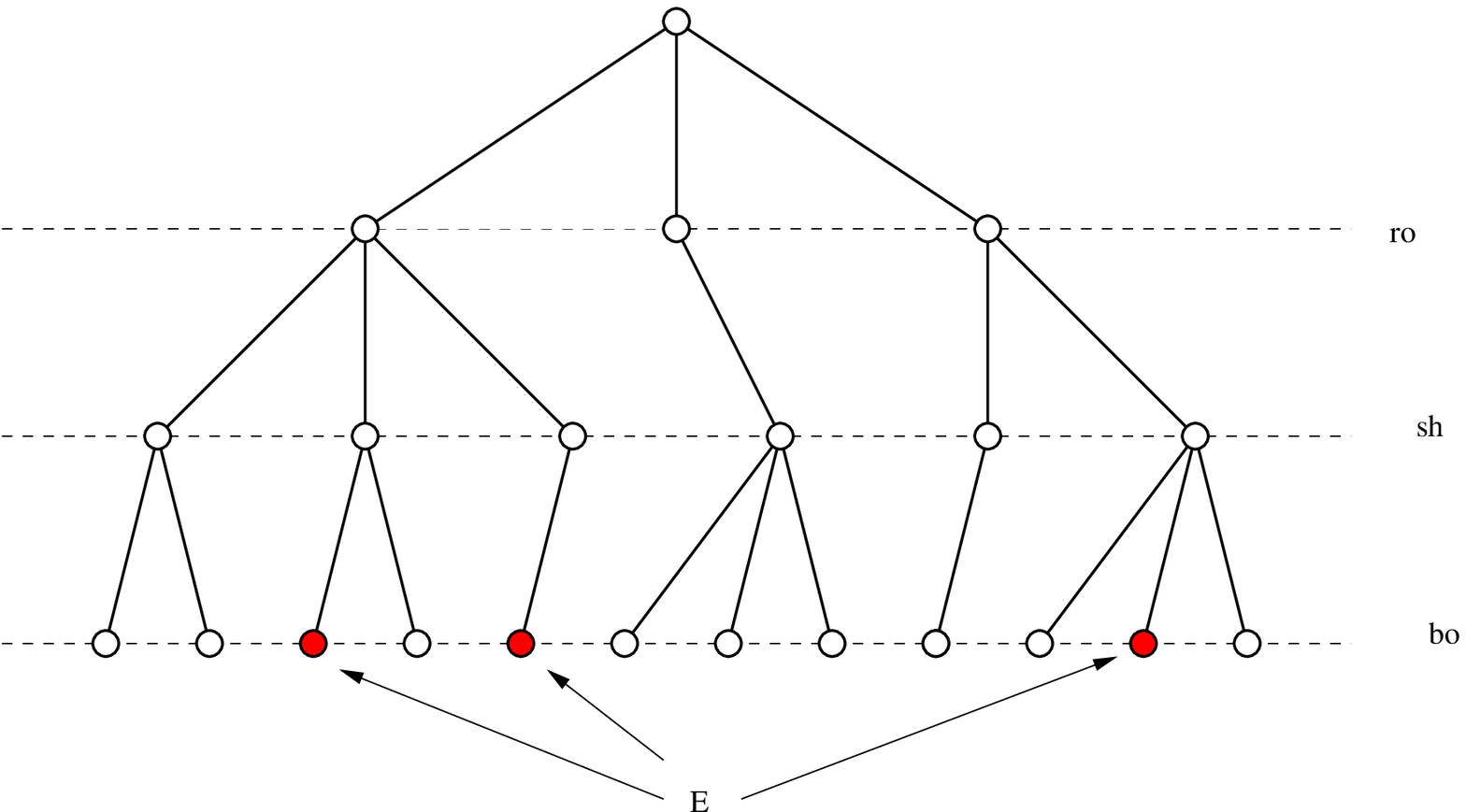}}
\end{center}
%\nopagebreak\vspace{.2cm} 
\centerline{{\bf Figure 2.}}
\vspace{3mm}

If  $\dep(T)=2$ we obtain the "library with several shelves" model considered in \cite{Bj},
as indicated in Figure 2. 
\noindent
This case was derived in \cite{Bj} via geometric considerations,
ultimately relying on Brown's theory of random walks on semigroups.
If one cares only about the library result, and not about random walks on complex hyperplane
arrangements, there is of course no need to mix in geometric considerations.
This note can be seen as a self-contained appendix to \cite{Bj} whose modest 
purpose is to fill in the details on how to obtain the
general  dynamic library model in the simplest and most direct way, avoiding geometry.

Another ``tree analogue''  of the classical linear random-to-front shuffle, different from
the one considered here,  has been studied in the literature. This is the random-to-root shuffle on
binary trees, see e.g.  \cite{AM,DF}.

\section{Shuffles on trees}

\newcommand{\Piord}{\Pi^{\text{ord}}}

%\section{Notation}
We begin by establishing notation.
For any finite set $A$, let
\beqarrs
S(A) &\defeq& \{\text{linear orderings of $A$}\} \\
\Pi(A) &\defeq& \{\text{partitions of $A$}\} \\
\Pi^{\text{ord}}(A) &\defeq& \{\text{ordered partitions of $A$}\} 
\eeqarrs
The   sets $\Pi(A)$ and $\Piord(A)$ are partially ordered by refinement,
meaning that $\al \le \be$ if and only if every block of the partition (or ordered partition) $\al$
is a union of blocks from $\be$.
Direct products (of sets, posets, \ldots)  are denoted by $\bigotimes$.

We consider rooted trees $T$ that are {\em pure}, meaning that all leaves are at
the same depth $d$. Let $V_j$ denote the set of nodes at depth $j$. So, $V_0 =
\{\text{root}\}$, $I\defeq \cup_{j=0}^{d-1} V_j=
\{\text{inner nodes}\}$, and $L\defeq V_d=\{\text{leaves}\}$.

\newcommand{\Cx}{C_{x}}
\newcommand{\Vx}{V_{\le x}}

\bddef
Let $E\sbseq L$.  A node $x\in T$ is {\em $E$-related} if some descendant of $x$
belongs to $E$.
%$\Vx \cap E\neq\emptyset$.
\eddef

For each inner node $x\in I$, let $\Cx$ denote the set of its children.
% and $\Vx$ the set of all descendants of $x$.
%The {\em degree}  at $x$ is $\deg (x)\defeq |\Cx|$.

\bddef
A {\em local ordering} of $T$ is a choice of linear order for the set of children $C_x$
at each inner node $x\in I$. Denote by $\mathcal{O} (T)$ the set of all local orderings of $T$.
Thus, $\OO (T) \cong \bigotimes_{x\in I} S(C_x)$.
\eddef

The subsets of $L$ act on $\OO(T)$ in the following way. 

\bddef \label{defmove}
Let $\pi = (\pi_x)_{x\in I}$
be a local ordering, and let $E\in 2^L$. Then
$E(\pi) = (E_x(\pi_x))_{x\in I}$, where
$E_x(\pi_x)$ is the linear ordering of $C_x$ in which the $E$-related elements come
first, in the order induced by $\pi_x$, followed by the remaining elements,
also in the induced order.
\eddef

The following figure shows a local ordering $\pi$ of a tree $T$, which 
coincides with left-to-right order in the planar drawing of $T$.
\vspace{.2cm}
%\btab{ccc}
\begin{center}
\psfrag{ro}{{rooms}}
\psfrag{sh}{{shelves}}
\psfrag{bo}{{books}}
\psfrag{E}{\large E}
\psfrag{2}{$2$}
\psfrag{3}{$3$}
\psfrag{1}{$1$}
\psfrag{x5}{\large {$y_5$}}
\psfrag{x6}{\large $y_6$}
\psfrag{x7}{\large {$y_7$}}
\psfrag{y1}{\large {$x_1 =x_5 =x_7$}}
\psfrag{y2}{\large $x_2 =x_3 =x_4 =x_6$}
\psfrag{R}{\large $\R$}
\resizebox{!}{50mm}{\includegraphics{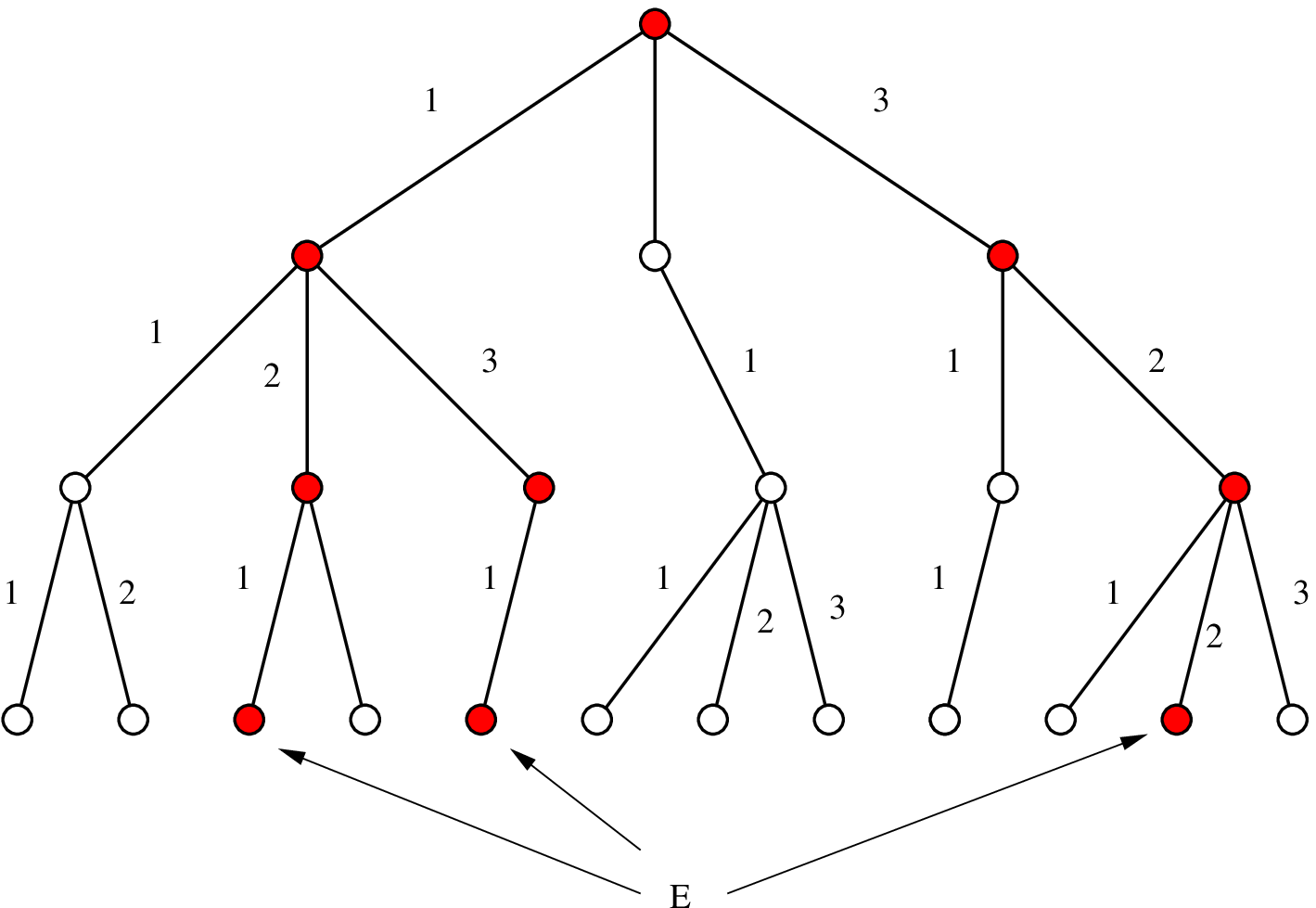}}
%{}{Archimedes}
\end{center}
%\nopagebreak\vspace{.2cm} 
\centerline{{\bf Figure 3.} }
\vspace{.2cm}
The indicated choice $E$ of leaves induces a move to the following local ordering $E(\pi)$.
The $E$-related nodes are shaded.% marked in red.
\vspace{.2cm}
%&$\Rightarrow$&
\begin{center}
\psfrag{ro}{{rooms}}
\psfrag{sh}{{shelves}}
\psfrag{bo}{{books}}
\psfrag{E}{\large E}
\psfrag{2}{$2$}
\psfrag{3}{$3$}
\psfrag{1}{$1$}
\psfrag{x5}{\large {$y_5$}}
\psfrag{x6}{\large $y_6$}
\psfrag{x7}{\large {$y_7$}}
\psfrag{y1}{\large {$x_1 =x_5 =x_7$}}
\psfrag{y2}{\large $x_2 =x_3 =x_4 =x_6$}
\psfrag{R}{\large $\R$}
\resizebox{!}{50mm}{\includegraphics{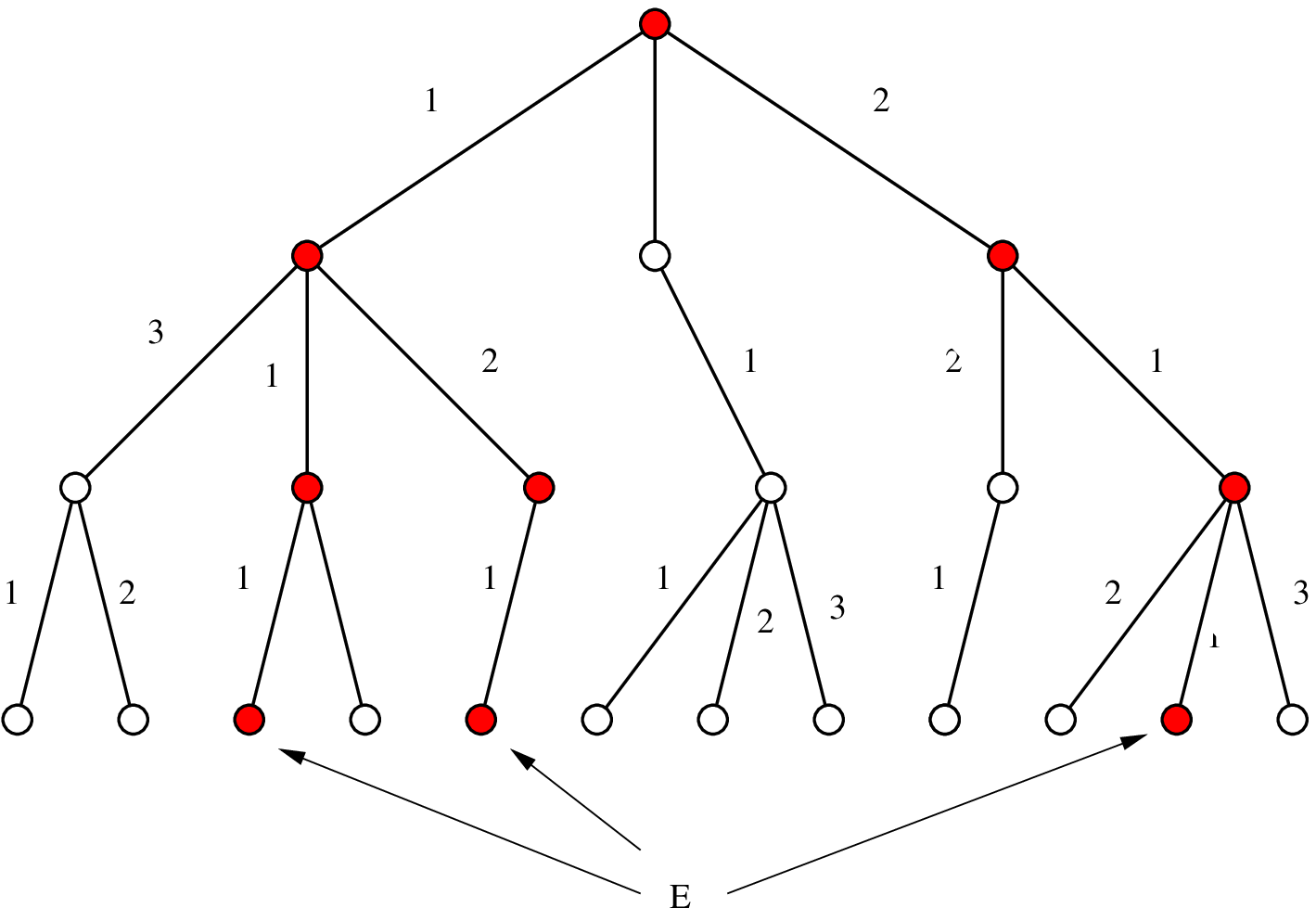}}
%{}{Archimedes}
\end{center}
%\etab
%\nopagebreak\vspace{.2cm} 
\centerline{{\bf Figure 4.} }
\vspace{2mm}

\bddef\label{walk}
Assume given a
probability distribution $(w_E)_{E\sbseq L}$ on $2^L$. This
determines a random walk on the set $\OO(T)$ as follows. 
If the walk is currently at the local ordering $\pi$, then choose a subset
$E\sbseq L$ with probability $w_E$ and move to $E(\pi)$.
\eddef

\newcommand{\PT}{\mathrm{Part}(T)}
Let $\PT \defeq \bigotimes_{x\in I}  \Pi(C_x)$. So, an element $\al\in\PT$ is a choice
of partition $\al_x$ of the set of children of $x$, for each inner node $x$.
The following special elements of $\PT$ are induced by subsets $E\sbseq L$.
For each $x\in I$ let $\al^E_x$ be the partition of
$C_x$ into two blocks, one block consisting of the $E$-related
elements and one of the remaining elements (one of these blocks may be
empty, in which case we forget it).

\bddef
Let $\al = (\al_x)_{x\in I} \in \PT$. 
A subset $E\sbseq L$ is {\em $\al$-compatible} if $\al_x$ is a refinement
of $\al^E_x$ for every $x\in I$.
\eddef
Notice that for every nontrivial $\al\in\PT$ there exists some $\al$-compatible
proper subset $E\sbseq L$.

\bthm \label{main}Let $T$ be a pure tree with leaves $L$.
Furthermore, let $\{w_E\}_{E\sbseq L}$ be a probability
distribution on $2^L$ and $P_w$ the transition matrix of the induced random walk on
local orderings of $T$:
$$P_w(\pi, \pi') =\sum_{E\st E(\pi)=\pi'} w_E
$$
for $\pi, \pi'\in \OO(T)$. Then,
\ben
\item[(i)] The matrix $P_w$ is diagonalizable.
%\vspace{1mm}
\item[(ii)] For each $\al =(\al_x)_{x\in I} \in \PT$ 
there is an eigenvalue
$$ \varepsilon_{\al} = \sum_{\mbox{E\st $E$ is $\al$-compatible}} w_E \, .
$$
%\vspace{2mm}
\item[(iii)] The multiplicity of the eigenvalue $\varepsilon_{\al}$ is 
$$m_{\al} = \prod_{x\in I} \prod_{B\in \al_x} (|B|-1)!
$$ 
\item[(iv)] These are all the eigenvalues of $P_w$.
\een
\ethm
\noindent
For clarity's sake, let us point out that
 $\ep_{\al}=\ep_{\be}$, for $\al\neq\be$, and $\ep_{\al}=0$
are possible.  

\bproof
As mentioned in the introduction, this is a special case of 
Brown's theory for walks on semigroups \cite{Bro1}, with which we now assume familiarity.

\newcommand{\PTord}{\mathrm{Part^{ord}}(T)}

Let $\PTord \defeq \bigotimes_{x\in I}  \Piord(C_x)$. So, an element $\be\in\PTord$ is a choice
of {\em ordered} partition $\be_x$ of the set of children of $x$, for each inner node $x$.
In particular, for each subset $E\sbseq L$ there is an element $\be^E \in \PTord$
whose component $\be^E_x$ at $x\in I$ is the two-block
ordered partition of $C_x$ whose first block consists of the $E$-related
elements of $C_x$, and second block of the remainder.
(If one of these blocks is empty we forget about it and let $\be_x^E$
have  only one block.)

Now, introduce the following probability distribution on $\PTord$:
\beq\mbox{Prob ($\be$) }  =
\bca w_E,  \mbox{ if }  \be=\be^E,\,\,  E\sbseq L \\
0,  \mbox{ \; for all other ordered partitions.}
\eca \eeq\label{prob}

Given this set-up, the proof consists of verifying each of the following claims for
$\PTord$, 
and then referring to %Theorem 1 on page 880 of  
\cite{Bro1}.

\ben
\item  $\PTord$ is an LRB (left regular band) semigroup with component-wise composition.
The composition in each factor $\Piord(A)$ has the following description.
If $X= \lang X_1, \dots, X_p\rang  $ and $Y=\lang Y_1, \dots, Y_q\rang  $ 
are ordered partitions of $A$,
%\; $X_i, Y_j \subseteq [n]$,
then $X\circ Y=\lang  X_i \cap Y_j\rang  $ with the blocks ordered by
the lexicographic order of the pairs of indices $(i,j)$.
%For instance,
%$$ \lang257\mid 3\mid 146\rang   \circ  \lang17\mid 25 \mid 346\rang   = \lang7\mid 25\mid
%3\mid 1\mid 46\rang . $$
\item
Its support lattice is $\PT$ and
support map
$$\supp: \PTord \rightarrow \PT,$$
whose component at each $x\in I$ is the map
$\Piord(C_x) \rarr\Pi(C_x)$ that sends an ordered partition of $C_x$
% $\lang   \ldots \rang  $
 to an unordered  partition 
 %$( \ldots )$
by forgetting the ordering of its blocks.
\item The maximal elements of $\PTord$ are the local orderings $\OO(T)$.
\item The steps of the 
semigroup random walk on $\OO(T)$, induced as in \cite{Bro1} by the probability
assignment (2.1),
%(\ref{prob}) 
are precisely the steps described in
Definition \ref{walk}.
\item For each $E\sbseq L$ and $\al\in\PT$:
$$ \supp(\be^E) \le \al \quad \Lrarr \quad E \mbox{ is $\al$-compatible}.
$$
\item The number of maximal elements of $\PTord$ above some $\be\in\PTord$ is
by Zaslavsky's theorem the sum of  M\"obius function absolute values $$\sum_{\al\ge\supp(\be)}
|\mu(\al, \ho)|$$
computed on the product partition lattice $\PT$. From this follows, via Brown's theory 
\cite{Bro1}, that
$$m_{\al}=|\mu(\al, \ho)|,$$
for all $\al\in\PT$.
By the product property of the \mob\ and its well-known explicit evaluation on
the partition lattice (see \cite{EC1}), this quantity equals 
$$|\mu(\al, \ho)| = \prod_{x\in I} \prod_{B\in \al_x} (|B|-1)!
$$ 
\een
\vspace{1mm}
In view of these facts 
the theorem is obtained by
specializing Theorem 1 on page 880 of \cite{Bro1}  to the semigroup $\PTord$.
\eproof

\section{Remarks}

\noindent 3.1. The random walk of Theorem \ref{main} has a unique stationary distribution $\pi$
 if and only if  $\{E\in 2^L \st w_E >   0\}$ is {\em separating}, meaning that
 for every inner node $x\in I$ and every pair of siblings $y,z \in C_x$, $y\neq z$, there is a subset
 $E\sbseq L$ with $w_E > 0$ for which one of $y$ and $z$ is $E$-related and the other is not.
 
 This follows from Theorem 2 of Brown and Diaconis \cite{BrDi}, using the fact that
 the random walk we consider can be realized as a walk on the complement of
 a product of real braid arrangements.
 Theorem 2 of \cite{BrDi} also gives  additional
information about the stationary distribution.
\vspace{3mm}

\noindent 3.2. One easily checks  that the subset $\{\be^E \st E\sbseq L \}$ 
generates the full semigroup $ {\mathrm{Part^{ord}}(T)}$, 
and that the set of its maximal
elements $\OO(T)$ is generated by \\ $\{\be^{\{e\}} \st e\in L \}$.
\vspace{3mm}

\noindent 3.3. Suppose that $w_E\neq 0$ only if $|E|=1$. Then Theorem \ref{main}
implies that the eigenvalues are indexed by $\bigotimes_{x\in I} 2^{C_x}$,
and that their multiplicities are products of derangement numbers, thus
generalizing the well-known result of
Donnelly, Kapoor-Reingold and Phatarfod for the Tsetlin library (the
$\dep(T)=1$ case); see the references for this given  in \cite{BHR, Bro1, BrDi}.


\begin{thebibliography}{}

\bibitem{AM} B. Allen and I.~Munro, {\em Self-organizing binary search trees},
J. Assoc. Comput. Mach. {\bf 25} (1978), 526--535.

\bibitem{BHR} P.~Bidigare, P.~Hanlon and D.~Rockmore,
{\em A combinatorial description of the spectrum for the Tsetlin library and its
generalization to hyperplane arrangements,} Duke Math J. {\bf 99} (1999), 135--174.

\bibitem{Bj} A.~Bj\"orner, {\em Random walks, arrangements, cell complexes, greedoids,
and self-organizing libraries,} in ``Building Bridges'' (eds. M. Gr\"otschel and G.  O. H. Katona), 
Bolyai Soc. Math. Studies {\bf 19} (2008),  Springer (Berlin) and Janos Bolyai Math. Soc. 
(Budapest), pp.165--203. 

\bibitem{Bro1} K.~S.~Brown, {\em Semigroups, rings and Markov chains,}
J. Theor. Probab. {\bf 13} (2000), 871--938.

\bibitem{BrDi} K.~S.~Brown and P. Diaconis, {\em Random walks and hyperplane arrangements},
Ann. Probab. {\bf 26} (1998), 1813--1854.

\bibitem{DF} R.~P.~Dobrow and J. A.~Fill, {\em On the Markov chain for the move-to-root rule 
for binary search trees}, Ann. Applied Probab. {\bf 5} (1995), 1--19.

\bibitem{FH} J. A.~Fill and L.~Holst, {\em On the distribution of search cost for the move-to-front rule},
Random Structures and Algorithms {\bf 8} (1996), 179--186.

\bibitem{EC1} R.~P.~Stanley, {\em Enumerative Combinatorics, Vol. 1}, Cambridge
Univ. Press, 1997.


\end{thebibliography}
\end{document}